\documentclass{article}
\usepackage[utf8]{inputenc}
\usepackage{amsmath}
\usepackage{amssymb}
\usepackage{amsthm}
\usepackage{enumerate}
\usepackage{mathrsfs}
\usepackage{scrextend}
\usepackage{mathtools}
\usepackage[hyperfootnotes=false]{hyperref}
\usepackage[multiple]{footmisc}
\usepackage[a4paper, total={6in, 9in}]{geometry}
\usepackage{scalerel,stackengine}
\stackMath
\newcommand\reallywidehat[1]{%
\savestack{\tmpbox}{\stretchto{%
  \scaleto{%
    \scalerel*[\widthof{\ensuremath{#1}}]{\kern-.6pt\bigwedge\kern-.6pt}%
    {\rule[-\textheight/2]{1ex}{\textheight}}
  }{\textheight}%
}{0.5ex}}%
\stackon[1pt]{#1}{\tmpbox}%
}
\theoremstyle{plain}
\newtheorem{theorem}{Theorem}

\newtheorem{proposition}[theorem]{Proposition}

\theoremstyle{definition}

\numberwithin{theorem}{section}

\title{A protrusive ordering of 5 points not witnessed by any finite multiset}
\author{Adrian Beker\footnote{University of Zagreb, Faculty of Science, Department of Mathematics, Zagreb,
Croatia.\\ Email: \nolinkurl{adrian.beker@math.hr}}}
\date{\today}

\begin{document}

\maketitle

\begin{abstract}
    Given a finite set of points $C \subseteq \mathbb{R}^d$, we say that an ordering of $C$ is \emph{protrusive} if every point lies outside the convex hull of the points preceding it. We give an example of a set $C$ of $5$ points in the Euclidean plane possessing a protrusive ordering that cannot be obtained by ranking the points of $C$ according to the sum of their distances to a finite multiset of points. This answers a question of Alon, Defant, Kravitz and Zhu.
\end{abstract}

\section{Introduction}

Given a point $x \in \mathbb{R}^d$ and a finite multiset $V = \{v_1, \ldots, v_k\}$ of points in $\mathbb{R}^d$, we let $D_V(x) = \sum_{i=1}^{k}\lVert x-v_i\rVert$ be the sum of the Euclidean distances of $x$ to the points of $V$. Given a finite set $C = \{c_1, \ldots, c_n\} \subseteq \mathbb{R}^d$ and a permutation $\sigma \in S_n$, we say that the ordering $(c_{\sigma(1)}, \ldots, c_{\sigma(n)})$ is \emph{witnessed} by $V$ if $D_V(c_{\sigma(1)}) < \ldots < D_V(c_{\sigma(n)})$. The set of orderings of $C$ witnessed by some finite multiset of points is denoted by $\Psi(C)$. 

A simple application of triangle inequality shows that any ordering $(c_{\sigma(1)}, \ldots, c_{\sigma(n)}) \in \Psi(C)$ must be \emph{protrusive}, that is to say $c_{\sigma(i+1)}$ must lie outside the convex hull of $c_{\sigma(1)}, \ldots, c_{\sigma(i)}$ for all $1 \leq i \leq n-1$. However, the converse does not hold, as was shown by Alon, Defant, Kravitz and Zhu in their recent preprint \hyperlink{adnz}{[1]}. They constructed a set of $6$ points $c_1, c_2, c_3, c_1', c_2', c_3'$ in $\mathbb{R}^2$ in convex position with the property that $D_{\{c_1,c_2,c_3\}}(v) \geq D_{\{c_1',c_2',c_3'\}}(v)$ for all $v \in \mathbb{R}^2$. In particular, the ordering $(c_1, c_2, c_3, c_1', c_2', c_3')$ is protrusive but is not in $\Psi(C)$. They also asked whether a set of $5$ points with the same property can be found (Question 7.5 in \hyperlink{adnz}{[1]}). The aim of this short note is to provide an affirmative answer to this question.

\begin{proposition}
\label{main result}
Let $C \subseteq \mathbb{R}^2$ consist of the five points
$$c_1 = (-1, -1), \quad\quad c_2 = (1, -1), \quad\quad c_3 = (0, 2),$$
$$c_1' = (0, 1.1), \quad\quad c_2' = (0, -1.1).$$
Then for any $v \in \mathbb{R}^2$, we have
\begin{equation}\label{main inequality}
    \frac{1}{3}D_{\{c_1,c_2,c_3\}}(v) \geq \frac{1}{2}D_{\{c_1',c_2'\}}(v).
\end{equation}
In particular, the ordering $(c_1, c_2, c_1', c_3, c_2')$ is protrusive but is not in $\Psi(C)$.
\end{proposition}

We remark that Proposition \ref{main result} is best possible in the sense that $C$ cannot be replaced by a set of fewer than five points in Euclidean space of any dimension nor by a set of any cardinality in $\mathbb{R}$ (this follows from Theorems 7.4 and 7.2 in \hyperlink{adnz}{[1]} respectively).

Both the statement and the proof of Proposition 1.1 bear significant similarity with those of Proposition 7.3 in \hyperlink{adnz}{[1]}. However, the important difference is that the set of points we construct is not in convex position. Hence, we do not have the freedom to consider any ordering of its points since not every ordering is protrusive. Nevertheless, this difficulty can be remedied by observing that a more general linear relation of the form $\sum_{i=1}^{n}\lambda_id_i \geq \sum_{j=1}^{m}\lambda_j'd_j'$ can be used to rule out more orderings of the numbers $d_1, \ldots, d_n, d_1', \ldots, d_m'$ than just those with all the $d_j'$ after all the $d_i$.

\section{Proof of Proposition \ref{main result}}

We follow the same strategy as in the proof of Proposition 7.3 in \hyperlink{adnz}{[1]}. We first show that (\ref{main inequality}) holds whenever $\lVert v\rVert$ is sufficiently large and we then check the remaining cases numerically. To begin, we record a slight generalisation of the first part of the argument from \hyperlink{adnz}{[1]}. Even though the proof follows very closely the argument in \hyperlink{adnz}{[1]}, we include it in full in order to make the note self-contained.

\begin{proposition}
\label{main result for large vectors}
Let $c_1, \ldots, c_n, c_1', \ldots, c_m' \in \mathbb{R}^2$ be such that $\sum_{i=1}^{n}c_i = 0 = \sum_{j=1}^{m}c_j'$. Let $R = \max_{i=1}^{n}\lVert c_i\rVert$, $R' = \max_{j=1}^{m}\lVert c_j'\rVert$ and suppose that $R,R' > 0$. Suppose that $\alpha > 1$ is such that the quadratic form
$$q(u) = \frac{1}{n}\sum_{i=1}^{n}\langle c_i,u\rangle^2 - \alpha \cdot \frac{1}{m}\sum_{j=1}^{m}\langle c_j',u\rangle^2$$
is positive semidefinite. Then
$$\frac{1}{n}\sum_{i=1}^{n}\lVert v-c_i\rVert \geq \frac{1}{m}\sum_{j=1}^{m}\lVert v-c_j'\rVert$$
holds for all $v \in \mathbb{R}^2$ such that $\lVert v\rVert > \max\left\{R,\frac{R+\alpha R'}{\alpha-1}\right\}$.
\end{proposition}
\noindent\textit{Proof.} Fix a vector $v \in \mathbb{R}^2$ such that $\lVert v\rVert > \max\left\{R,\frac{R+\alpha R'}{\alpha-1}\right\}$, so in particular $\lVert v\rVert > R'$. Choose a unit vector $u \in \mathbb{R}^2$ perpendicular to $v$ and let $\pi \colon \mathbb{R}^2 \to \mathbb{R}^2$ denote orthogonal projection onto the line spanned by $v$. Given $x \geq 0$ and $a > 0$, we observe the elementary estimates
$$a + \frac{x}{2\sqrt{a^2+x}} \leq \sqrt{a^2+x} \leq a + \frac{x}{2a},$$
which follow, say, by writing
$$\sqrt{a^2+x}-a = \frac{\left(\sqrt{a^2+x}-a\right)\left(\sqrt{a^2+x}+a\right)}{\sqrt{a^2+x}+a} = \frac{x}{\sqrt{a^2+x}+a}.$$
Given any $c \in \mathbb{R}^2$ such that $\pi(c) \neq v$, on applying these estimates with $a = \lVert v-\pi(c)\rVert$, $x = \lVert \pi(c)-c\rVert^2$ and using Pythagoras' theorem, we find that
\begin{equation}\label{auxiliary estimate}
    \lVert v-\pi(c)\rVert + \frac{\lVert \pi(c)-c\rVert^2}{2\lVert v-c\rVert} \leq \lVert v-c\rVert \leq \lVert v-\pi(c)\rVert + \frac{\lVert \pi(c)-c\rVert^2}{2\lVert v-\pi(c)\rVert}.
\end{equation}
Since $\lVert c_i\rVert \leq \lVert v\rVert$, the Cauchy-Schwarz inequality implies that $|\langle c_i,v\rangle| \leq \lVert v\rVert^2$. Hence, we have
$$\frac{1}{n}\sum_{i=1}^{n}\lVert v-\pi(c_i)\rVert = \frac{1}{n}\sum_{i=1}^{n}\left\lVert v-\frac{\langle c_i,v\rangle}{\lVert v\rVert^2}v\right\rVert = \frac{1}{n}\sum_{i=1}^{n}\left(1-\frac{\langle c_i,v\rangle}{\lVert v\rVert^2}\right)\lVert v\rVert = \lVert v\rVert - \frac{1}{n}\left\langle \sum_{i=1}^{n}c_i, \frac{v}{\lVert v\rVert}\right\rangle = \lVert v\rVert.$$
Similarly, it follows that
$$\frac{1}{m}\sum_{j=1}^{m}\lVert v-\pi(c_j')\rVert = \lVert v\rVert.$$
By combining this with the estimate (\ref{auxiliary estimate}) and using the triangle inequality, we obtain that
$$\frac{1}{n}\sum_{i=1}^{n}\lVert v-c_i\rVert \geq \frac{1}{n}\sum_{i=1}^{n}\lVert v-\pi(c_i)\rVert + \frac{1}{n}\sum_{i=1}^{n}\frac{\lVert \pi(c_i)-c_i\Vert^2}{2\lVert v-c_i\rVert} \geq \lVert v\rVert + \frac{1}{n}\sum_{i=1}^{n}\frac{\langle c_i,u\rangle^2}{2(\lVert v\rVert + R)},$$
$$\frac{1}{m}\sum_{j=1}^{m}\lVert v-c_j'\rVert \leq \frac{1}{m}\sum_{j=1}^{m}\lVert v-\pi(c_j')\rVert + \frac{1}{m}\sum_{j=1}^{m}\frac{\lVert \pi(c_j')-c_j'\rVert^2}{2\lVert v-\pi(c_j')\rVert} \leq \lVert v\rVert + \frac{1}{m}\sum_{j=1}^{m}\frac{\langle c_j',u\rangle^2}{2(\lVert v\rVert - R')}.$$
On subtracting the obtained inequalities, it follows that
$$\frac{1}{n}\sum_{i=1}^{n}\lVert v-c_i\rVert - \frac{1}{m}\sum_{j=1}^{m}\lVert v-c_j'\rVert \geq \frac{1}{n}\sum_{i=1}^{n}\frac{\langle c_i,u\rangle^2}{2(\lVert v\rVert + R)} - \frac{1}{m}\sum_{j=1}^{m}\frac{\langle c_j',u\rangle^2}{2(\lVert v\rVert - R')}.$$
Since $\frac{\lVert v\rVert + R}{\lVert v\rVert-R'} \leq \alpha$, the right-hand side can be bounded from below as follows:
$$\frac{1}{2(\lVert v\rVert+R)}\left(\frac{1}{n}\sum_{i=1}^{n}\langle c_i,u\rangle^2 - \frac{\lVert v\rVert + R}{\lVert v\rVert-R'} \cdot \frac{1}{m}\sum_{j=1}^{m}\langle c_j',u\rangle^2\right) \geq \frac{q(u)}{2(\lVert v\rVert+R)} \geq 0,$$
so we are done. $\qed$
\\\\
\noindent\textbf{Remark.} By making a clever choice of the $c_i$ and the $c_j'$, both terms in the analogue of the quadratic form $q$ that appears in \hyperlink{adnz}{[1]} turn out to be multiples of the square of the standard Euclidean norm. As a consequence, the corresponding quadratic form is identically equal to zero. Our contribution lies in relaxing this condition by noticing that it is enough for $q$ to be positive semidefinite, which is a condition that is easily checkable in practice.
\bigskip

We now turn to the proof of the main result, which is computer assisted. The relevant code can be found on the following link:

\begin{center}
    \url{https://github.com/adrianbeker/protrusive-unwitnessed}.
\end{center}

\noindent\textit{Proof of Proposition \ref{main result}.} We see that $c_1+c_2+c_3 = 0 = c_1'+c_2'$ and $R \vcentcolon= \max\{\lVert c_1\rVert,\lVert c_2\rVert,\lVert c_3\rVert\} = 2$, $R' \vcentcolon= \max\{\lVert c_1'\rVert, \lVert c_2'\rVert\} = 1.1$. Furthermore, letting $\alpha \vcentcolon= \frac{2}{1.1^2} > 1$, we have
\begin{align*}
    q(x,y) \vcentcolon= \frac{(-x-y)^2 + (x-y)^2 + (2y)^2}{3} - \alpha \cdot \frac{(1.1y)^2 + (-1.1y)^2}{2} = \frac{2x^2 + 6y^2}{3} - 1.1^2\alpha y^2 = \frac{2}{3}x^2,
\end{align*}
which is a positive semidefinite quadratic form. Since $\frac{R+\alpha R'}{\alpha-1} \approx 5.85 < 6$, Proposition \ref{main result for large vectors} implies that (\ref{main inequality}) holds whenever $\lVert v\rVert \geq 6$. It thus suffices to check that (\ref{main inequality}) holds whenever $\lVert v\rVert \leq 6$. To this end, we introduce the function
$$f \colon \mathbb{R}^2 \to \mathbb{R}, \quad v \mapsto 2\sum_{i=1}^{3}\lVert v-100c_i\rVert - 3\sum_{j=1}^{2}\lVert v-100c_j'\rVert$$
and observe that it suffices to check that $f$ is non-negative on $[-600,600]^2$. If we consider instead the integer-valued function
$$\widetilde{f} \colon \mathbb{R}^2 \to \mathbb{Z}, \quad v \mapsto 2\sum_{i=1}^{3}\left\lfloor\lVert v-100c_i\rVert\right\rfloor - 3\sum_{j=1}^{2}\left\lceil\lVert v-100c_j'\rVert\right\rceil,$$
a simple computer program (avoiding the use of floating-point arithmetic altogether) can be used to verify that the minimum of $\widetilde{f}$ on $([-600,600] \cap \mathbb{Z})^2$ is $12$. But $f$ is $12$-Lipschitz and $([-600,600] \cap \mathbb{Z})^2$ is a $\frac{\sqrt{2}}{2}$-net for $[-600,600]^2$, so since $f \geq \widetilde{f}$, it follows that
$$f(v) \geq 12 - 12 \cdot \frac{\sqrt{2}}{2} > 0$$
for all $v \in [-600,600]^2$, as desired. Finally, since $c_1'$ lies strictly inside the triangle formed by $c_1, c_2, c_3$ and $c_2'$ lies strictly outside it, it follows that the ordering $(c_1,c_2,c_1',c_3,c_2')$ is protrusive. Moreover, if $V = \{v_1,\ldots,v_k\}$ is a finite multiset of points in $\mathbb{R}^2$, then on summing the inequality (\ref{main inequality}) with $v_l$ in place of $v$ over all $l \in [k]$ and interchanging the order of summation, we obtain that
$$\frac{1}{3}\sum_{i=1}^{3}D_V(c_i) \geq \frac{1}{2}\sum_{j=1}^{2}D_V(c_j').$$
But if $V$ were to witness $(c_1,c_2,c_1',c_3,c_2')$, then we would have
$$\frac{1}{3}\sum_{i=1}^{3}D_V(c_i) < \frac{1}{3}(2D_V(c_1')+D_V(c_2')) < \frac{1}{2}\sum_{j=1}^{2}D_V(c_j'),$$
which is absurd. Therefore, this ordering does not belong to $\Psi(C)$, as desired. $\qed$
\bigskip

\noindent\textbf{Acknowledgements.} The author would like to thank Rudi Mrazović and Daniel G. Zhu for helpful comments on a draft of the note.

\end{document}